\documentclass[12pt,twoside]{article}

\usepackage{latexsym}
\usepackage[applemac]{inputenc}
\usepackage{fancyhdr}
\usepackage{amsmath, amsthm, amssymb, amsfonts, enumerate}
\usepackage[applemac]{inputenc} 
\usepackage[colorlinks=true,linkcolor=blue,urlcolor=blue]{hyperref}

\newtheorem{thm}{Theorem}[section]
\newtheorem{cor}[thm]{Corollary}
\newtheorem{lem}[thm]{Lemma}
\newtheorem{pro}[thm]{Proposition}
\newtheorem{defi}[thm]{Definition}

\def \Proof{\noindent{\bf Proof:} \ }

\def \conv{ {\rm conv} }

\def \conv{ {\rm conv} }

\def \lim{ {\rm lim} }

\def \Proof{\noindent{\bf Proof:} \ }

\def \conv{ {\rm conv} }

\def \conv{ {\rm conv} }

\def \conv{ {\rm conv\,} }

\def \R{I\!\!R}
\def \E{I\!\!E}
\def \P{ {\bf P} }

\def \N{I\!\!N}

\def \R{I\!\!R}
\def \E{I\!\!E}
\def \P{ {\bf P} }

\def \N{I\!\!N}
 
 \def \1{{\bf 1}}

\def \dis{\displaystyle}

\def \conv{ {\rm conv} }

\def \Val{ {\rm Val} }

 \title{A tutorial on Zero-sum Stochastic Games}

 \author{J\'er\^ome Renault\thanks{Toulouse School of Economics, University of Toulouse Capitole, Toulouse, France. This tutorial  has been originally prepared for  the program ``Stochastic Methods in Game Theory" which took place from November 16 to  December 25, 2015, at the National University of Singapore.  The lectures 
were  completed while the author was visiting the Institute for Mathematical Sciences, National University of Singapore in 2015. Repeated Games with incomplete information  and a few open problems have also been presented during the tutorial, but do not appear in the present version of this document. }}
\date{\today}

\begin{document}

 \maketitle

\begin{abstract} Zero-sum stochastic games generalize the notion of Markov Decision Processes (i.e. controlled Markov chains, or stochastic dynamic programming) to the 2-player competitive case : two players jointly control the evolution of a state variable, and have opposite interests.  These notes constitute a  short mathematical introduction to the theory  of such  games. Section 1 presents the basic model  with finitely many states and actions.  We give proofs of the standard results concerning : the existence and formulas for the values of the $n$-stage games, of the $\lambda$-discounted games, the  convergence of these values when $\lambda$ goes to 0 (algebraic approach) and when $n$ goes to $+\infty$, an important example called``The Big Match" and the existence of the uniform value. Section 2 presents a short and subjective selection of related and more recent results  : 1-player games (MDP)   and the compact non expansive case,  a simple compact continuous stochastic game with no asymptotic value, and  the general equivalence between the uniform convergence of $(v_n)_n$ and $(v_\lambda)_\lambda$.

More references on the topic can be found for instance in the books by Mertens-Sorin-Zamir (1994, reedited in 2015), Sorin (2002), Neyman-Sorin (2003)  or the chapters by Vieille (2002) and Laraki-Sorin (2014).\end{abstract}

%
%
%
%

%
%
%
%
%
%
%
%
%
%
%
%
%
 
 
 \section{The basic model}
 
 \subsection{Description}
 
 Zero-sum games are 2-player games where the sum of the payoffs of the players is 0, they are games of pure competition between the players.
 Zero-sum stochastic games are   dynamic zero-sum  games played in discrete time.  The basic model is due to Shapley (1953),  and is given by a set of states $K$ with an initial state $k_1$,  a set of actions $I$ for player 1, a set of actions $J$ for player 2, a payoff function $g:K \times I \times J \longrightarrow \R$, and a transition mapping $q$ from $K \times I \times J$ to the simplex $\Delta(K)$ of probability distributions over $K$. In the basic  model, $K$, $I$ and $J$ are assumed  to be non empty {\it  finite} sets.\\

 The progress of the game is the following:
 
 - The initial state is $k_1$, known to the players. At stage 1, player 1 and player 2  simultaneously choose  $i_1 \in I$ and $j_1 \in J$. Then P1's payoff is  $g(k_1,i_1,j_1)$ and P2 's payoff is   $-g(k_1,i_1,j_1)$, the actions  $i_1$ et $j_1$ are publicly announced,  and the play proceeds to stage 2.

- A stage $t\geq 2$, the state $k_t$ is selected according to the distribution  $q(k_{t-1},i_{t-1},j_{t-1})$, and is announced to both players. Player 1 and player 2  then simultaneously choose  $i_t \in I$ et $j_t \in J$.  P1's payoff is  $g(k_t,i_t,j_t)$ and P2's payoff is  $-g(k_t,i_t,j_t)$, the actions $i_t$ et $j_t$ are announced, and the play proceeds to stage   $t+1$.\\

\noindent{\bf Notations and  vocabulary.} We denote by $q(k'|k,i,j)$ the probability that the state of stage $t+1$ is $k'$ if the state of stage $t$ is $k$ and  $i$ and $j$ are played at that stage. A state  $k$ is absorbing if  $q(k|k,i,j)=1$ for all  $(i,j)$ in  $I \times J$ (when $k$ is reached, the play stays there forever). A stochastic game is absorbing if it has a unique non absorbing state.

A play is a sequence  $(k_1,i_1,j_1,k_2,i_2,j_2,....,k_t,i_t,j_t,...)$  taking values in $K \times I \times J$.\\
A history of the game is a finite sequence $(k_1,i_1,j_1,....,k_{t-1},i_{t-t},j_{t-1},k_t)$ in $(K\times I \times J)^{t-1}\times K$  for some positive integer $t$,  representing the information available to the players before they play at stage $t$.

A behavior strategy, or simply a strategy  of player 1 (resp. player 2), associates to every history a mixed action in $\Delta(I)$ (resp. $\Delta(J)$) to be played in case  this history occurs. 
A strategy of a player is said to be pure if it associates to each history a Dirac measure, that is an element of $I$ for player 1 and an element of $J$ for player 2.

A strategy is said to be Markov  if for any stage $t$, the mixed action prescribed at stage $t$ only depends on the current state $k_t$ (and not on past states or past actions).
A stationary strategy is a Markov  strategy such that the mixed action prescribed after any history only depends on the current state (and not on the stage number).\\

 As usual,  in all the examples  player 1 chooses the row, player 2 the column and we indicate the payoff   of player 1 (which  is minus the payoff of player 2). An absorbing state will be denoted with a *. For instance, 3* represents an absorbing state where the payoff to player 1 is 3, whatever the actions played.\\

\noindent{\bf Example 1:}
\begin{center}
\hspace{2,5cm} $L$ \hspace{0,2cm} $R$ \hspace{1,5cm}  
\vspace{0,2cm}

$\dis
\begin{array}{c}
T\\
B
\end{array}
\left( \begin{array}{cc}
0  & 1^*\\
1^* & 0^* \\
\end{array}
\right)$

\end{center}

There is a unique non absorbing state which is the initial state. Actions are $T$ and $B$ for player 1, $L$ and $R$ for player 2. If at the first stage the action profile played is $(T,L)$ then the stage payoff is 0 and the play goes to the next stage without changing state. If at the first stage the action profile played is $(T,R)$ or $(B,L)$, the play reaches an absorbing state where at each subsequent stage, whatever the actions played the payoff of player 1 will be 1.  If at the first stage the action profile played is $(B,R)$, the play reaches an absorbing state where at each subsequent stage, whatever the actions played the payoff of player 1 will be 0. \\

 \noindent {\bf  Example 2:} A one-player game ($J$ is a singleton), with determinisitic transitions and actions Black and Blue for Player 1. The payoffs are either 1 or 0 in each case.
 
 \vspace{0,5cm}
 
  \unitlength 0,7mm
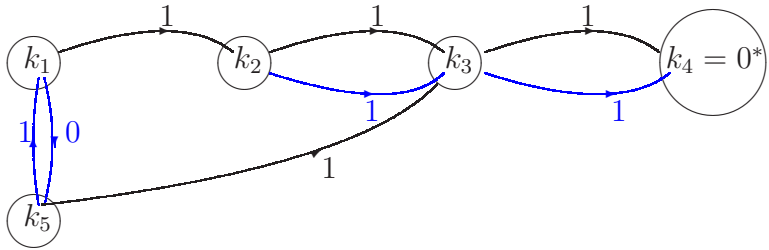
 \begin{center}
\begin{picture}(200,35)
 \put(00,0){\circle{10}}
   \put(00,30){\circle{10}}
    \put(40,30){\circle{10}}
       \put(80,30){\circle{10}}
    \put(129,30){\circle{20}}

    \put(-2,-1){$k_5$}
    \put(-2,29){$k_1$}
    \put(38,29){$k_2$}
    \put(78,29){$k_3$}
    \put(120,29){$k_4=0^*$}
    
  \textcolor{blue} { \qbezier(2,27)(5,15)(2,3)}
  \textcolor{blue} { \put(2,15){\vector(0,-1){1}}}
    \textcolor{blue} {   \put(2,15){$0$}}

  \textcolor{blue} {     \qbezier(-5,27)(-7,15)(-5,3)}
  \textcolor{blue} { \put(-8,15){\vector(0,1){1}}}
   \textcolor{blue} {   \put(-13,15){$1$}}

  { \qbezier(-5,32)(18,40)(28,32)}
 { \put(15,36){\vector(1,0){1}}}
    {   \put(14,37){$1$}}
    
      { \qbezier(35,32)(58,40)(68,32)}
 { \put(55,36){\vector(1,0){1}}}
    {   \put(54,37){$1$}}
    
  \textcolor{blue} {          { \qbezier(33,28)(58,20)(66,28)}}
  \textcolor{blue} {    { \put(50,24){\vector(1,0){1}}}}
  \textcolor{blue} {       {   \put(47,19){$1$}}}
    
   { \qbezier(70,32)(93,40)(103,32)}
 { \put(90,36){\vector(1,0){1}}}
    {   \put(88,37){$1$}}

    \textcolor{blue} {         { \qbezier(68,28)(93,20)(103,28)}}
 \textcolor{blue} {    { \put(90,24){\vector(1,0){1}}}}
  \textcolor{blue} {      {   \put(88,19){$1$}}}
  
  
%
  
   {         { \qbezier(-20,3)(40,10)(55,26)}}
 {    { \put(32,13){\vector(2,1){1}}}}
  {      {   \put(33,8){$1$}}}

\end{picture}
 
 \end{center}

  \vspace{0,5cm}
  
 \noindent {\bf  Example 3:}    The ``Big Match" 
 
 \centerline{$ \left(\begin{array}{cc}
 1^* & 0^* \\
0 & 1 \\
\end{array}
\right) $}

\vspace{1cm}

We denote by ${\Sigma}$ and ${\cal T}$ the sets of strategies of player 1 and 2, respectively.   A  couple of strategies in $\Sigma \times {\cal T}$ naturally\footnote{just as   tossing  a coin at every stage induces a probability distribution over sequences of Heads and Tails.} induces a probability distribution $\P_{k_1,\sigma, \tau}$ over the set of plays $\Omega=(K \times I \times J)^\infty$, endowed with the product $\sigma$-algebra. We will denote the expectation with respect to $\P_{k_1,\sigma, \tau}$ by $\E_{k_1,\sigma, \tau}$.\\

\noindent{\bf Remark:} A mixed strategy of a player is a probability distribution over his set of pure strategies (endowed with the product $\sigma$-algebra). By Kuhn's theorem (Aumann, 1964), one can show that mixed strategies and behavior strategies are equivalent, in the following strong sense : for any behavior  strategy $\sigma$ of player 1 there exists a mixed strategy $\sigma'$ of this player such that, for any pure (or mixed, or behavior) strategy of player 2, $(\sigma, \tau)$ and $(\sigma', \tau)$ induce the same probabilities over plays. And vice-versa by exchanging the words ``mixed" and  ``behavior" in the last sentence. Idem by exchanging the roles of player 1 and player 2   above. \\
  
\subsection{The  $n$-stage game and the $\lambda$-discounted game}

\begin{defi} \label{def1}Given a positive integer $n$, the $n$-stage game with initial state $k_1$ is the zero-sum game $\Gamma_n(k_1)$ with strategy spaces $\Sigma$ and ${\cal T}$ payoff function:
$$\forall (\sigma, \tau)\in  \Sigma \times {\cal T}, \;\; \gamma_n^{k_1}(\sigma,\tau)=\E_{k_1, \sigma, \tau} \;\left( \frac{1}{n} \sum_{t=1}^n g(k_t,i_t,j_t)\right).$$\end{defi} 

Because only finitely many stages matter here,  $\Gamma_n(k_1)$ can be equivalently seen as a finite zero-sum game played with mixed strategies. Hence it has a value denoted by $v_n(k_1)=\max_{\sigma\in \Sigma}\min_{\tau \in {\cal T}} \gamma_n^{k_1}(\sigma,\tau)$ $=\min_{\tau \in {\cal T}}\max_{\sigma\in \Sigma} \gamma_n^{k_1}(\sigma,\tau).$ For convenience we write $v_0(k)=0$ for each $k$. 

\begin{defi}  \label{def2} Given a discount factor $\lambda$ in $(0,1]$, the $\lambda$-discounted game  with initial state $k_1$ is the zero-sum game $\Gamma_\lambda(k_1)$ with strategy spaces $\Sigma$ and ${\cal T}$ payoff function:
$$\forall (\sigma, \tau)\in  \Sigma \times {\cal T}, \;\; \gamma_\lambda^{k_1}(\sigma,\tau)=\E_{k_1, \sigma, \tau} \;\left( \lambda \sum_{t=1}^{\infty}(1-\lambda)^{t-1}  g(k_t,i_t,j_t)\right).$$\end{defi} 

By a variant of Sion theorem, it has a value denoted by $v_\lambda(k_1)$. In the economic literature $\delta=1-\lambda=\frac{1}{1+r}$ is called the discount rate, $r$ being called the interest rate.

\begin{pro} \label{pro1} $v_n$ and $v_\lambda$ are characterized by the following Shapley equations.\\

1) For $n\geq 0$ and  $k$ dans $K$:
$$(n+1)\;  v_{n+1} (k)=\Val_{\Delta(I)\times \Delta(J)}  \left(g(k,i,j)+ \sum_{k' \in K} q(k'|k,i,j)\;  n \;  v_n(k')\right).$$
And in any $n$-stage game, players have Markov optimal  strategies.\\

2) For $\lambda$ in $(0,1]$ and $k$ in $K$:
$$  v_{\lambda} (k)=\Val_{\Delta(I)\times \Delta(J)}  \left(\lambda \; g(k,i,j)+(1-\lambda)  \sum_{k' \in K} q(k'|k,i,j)\;   \;  v_{\lambda}(k')\right).$$
And in any  $\lambda$-discounted game, players have stationary optimal strategies.

\end{pro}

 \noindent{\Proof} The proof is standard. For 1), fix $n$ and $k$ and denote by $v$ the value of the matrix game $ \left(g(k,i,j)+ \sum_{k' \in K} q(k'|k,i,j)\;  n \;  v_n(k')\right)_{i,j}$. In the game with $n+1$ stages and initial state $k$, player 1 can play at stage 1 an optimal strategy in this matrix game, then from stage 2 on an optimal strategy in the remaining $n$-stage stochastic game. By doing so, player 1 guarantees $v$ in  $\Gamma_{n+1}(k)$, so $v_{n+1}(k)\geq v$. Proceeding similarly with player 2 gives $v_{n+1}(k)= v$.
 
 The proof of 2) is similar. Notice that by the contracting fixed point theorem, for fixed $\lambda$  the vector $(v_\lambda(k))_{k \in K}$ is uniquely characterized by the Shapley equations.\\

It is easy to compute $v_n$ and $v_\lambda$ in the previous examples (in absorbing games, we simply  write $v_n$ and $v_\lambda$ for the values of the stochastic game where the initial state is the non absorbing state)

Example 1: $v_1=\frac{1}{2}$, $v_{n+1}=\frac{1}{2- {\frac{n}{n+1}v_n}}$ for $n\geq 1$, and $v_\lambda=\frac{1}{1+\sqrt{\lambda}}$ for each $\lambda$.

Example 2: For $\lambda$ small enough, $v_\lambda(k_1)=\frac{1}{2-\lambda}$ and it is optimal in the $\lambda$-discounted game to alternate between states $k_1$ and $k_5$. For $n\geq 0$, $(2n+3)v_{2n+3}=(2n+4)v_{2n+4}={n+3}$ (first alternate between $k_1$ and $k_5$, then go to $k_2$ 3 or 4 stages before the end).

Example 3 (The Big Match): $v_n=v_\lambda=1/2$ for all $n$ and $\lambda$.\\

The Shapley operator is defined as the mapping which associates to each $v$ in $\R^K$ the vector $\Psi(v)$ in $\R^K$ such that for each $k$,
$${\Psi(v)}^k= \Val_{\Delta(I)\times \Delta(J)}  \left(g(k,i,j)+ \sum_{k' \in K} q(k'|k,i,j)\;    v^{k'}\right).$$
$\Psi$ is non expansive for the sup-norm $\|v\|=\sup_{k\in K} |v^k|$ on $\R^K$, and the Shapley equations can be rephrased as:

$$\forall n\geq 1,\, nv_{n}= \Psi((n-1) v_{n-1})=\Psi^n(0),$$

$$\forall \lambda\in (0,1], v_\lambda=\lambda \Psi\left(\frac{1-\lambda}{\lambda}v_\lambda\right).$$

\subsection{Limit values - The algebraic approach}

We are  interested here in the limit values when the players become more and more patient, i.e. in the existence of the limits of $v_n$, when $n$ goes to infinity, and of $v_\lambda$, when $\lambda$ goes to 0.

It is always interesting to study first the 1-player case.

\subsubsection{1-player case: Markov Decision Process}

We assume here that player 2 does not exist, that is $J$ is a singleton.
For any  $\lambda>0$, player 1 has an optimal  stationary strategy in the $\lambda$-discounted game. Moreover since the matrix games appearing in 2) of the Shapley equations only have one column, this stationary optimal strategy  can   be taken to be pure. So we just have to consider strategies given by a mapping $f:K\longrightarrow I$, with the interpretation that player 1 plays $f(k)$ whenever the current state is $k$.

The $\lambda$-discounted payoff when $f$ is played and the initial state is $k$ satisfies:

$$\gamma_\lambda^k(f)= \lambda g(k,f(k))+ (1-\lambda) \sum_{k'\in K} q(k'|k,f(k)) \gamma_\lambda^{k'}(f).$$

\noindent Consider the vector $v=(\gamma_\lambda^k(f))_k$. The above equations can be written in matrix form: $(I-(1-\lambda)A)v= \lambda \alpha$, where $I$ is the identity matrix, $A=(q(k'|k,f(k))_{k,k'}$ is a stochastic matrix independent of $\lambda$,  and $\alpha=(g(k,f(k)))_k$ is a fixed vector. $(I-(1-\lambda)A)$ being invertible, we know that its inverse has coefficients which are rational fractions of its coefficients. Consequently, we obtain that:

For a given pure stationary strategy $f$, the payoff $\gamma_\lambda^k(f)$ is a rational fraction of $\lambda$.\\

Now we have finitely many such strategies to consider, and a given $f$ is optimal in the $\lambda$-discounted game with initial state $k$ if and only if:
 $\gamma_\lambda^k(f)\geq  \gamma_\lambda^k(f')$ for all $f'$. Because a non-zero polynomial only has finitely many roots, we obtain that for $\lambda$ small enough, the same pure optimal strategy $f$ has to be optimal in any discounted game. And clearly $f$ can be taken to be optimal whatever the initial state is.
 
 \begin{thm} (Blackwell, 1962)
 In the 1-player case, there exists $\lambda_0>0$ and a pure stationary optimal strategy $f$ which is optimal in any $\lambda$-discounted game with $\lambda\leq \lambda_0$.
 For $\lambda\leq \lambda_0$ and $k$ in $K$, the value $v_\lambda(k)$ is a bounded rational fraction of $\lambda$, hence converges when $\lambda$ goes to 0.
 \end{thm}
 
 In example 2, $f$ is the strategy which alternates forever between $k_1$ and $k_5$. There exists no strategy which is optimal in all $n$-stage games  with $n$ sufficiently large.

\subsubsection{Stochastic games: The algebraic approach}

Back to the 2-player case, we know that in each discounted game the players have stationary optimal strategies.
The following approach\footnote{M. Oliu-Barton (2014) provided a  proof of the convergence of $v_\lambda$ using elementary tools.}   is due to Bewley and Kohlberg (1976). Consider the following set:

\begin{eqnarray*}
A = & \{(\lambda, x_\lambda, y_\lambda, w_\lambda) \in (0,1]\times (\R^I)^K \times (\R^J)^K \times \R^K, \forall k\in K,\\
\; &  \; x_\lambda(k), y_{\lambda}(k) \; {\rm stationary \; optimal \;in\; } \Gamma_\lambda(k),\; w_\lambda(k)=v_{\lambda}(k)\}.
 \end{eqnarray*}
 
 $A$ can be written with finitely many polynomial inequalities:
 
 $$\forall i,j,k,\; \sum_{i} x_\lambda^i(k)=1, x_\lambda^i(k) \geq 0, \sum_j y_\lambda^j(k)=1, y_\lambda^j(k)\geq 0,$$
 
 $$\forall j,k,\; \sum_{i\in I} x_\lambda^i(k) (\lambda g(k,i,j)+(1-\lambda) \sum_{k'} q(k'|k,i,j) w_{\lambda}(k')) \geq w_{\lambda}(k),$$
 
  $$\forall i,k,\; \sum_{j\in J} y_\lambda^j(k) (\lambda g(k,i,j)+(1-\lambda) \sum_{k'} q(k'|k,i,j) w_{\lambda}(k')) \leq w_{\lambda}(k).$$

In particular, the set $A$ is    semi-algebraic\footnote{A subset of an Euclidean space is semi-algebraic if it can be written  a finite union of sets, each of these sets  being defined as the conjunction of finitely many weak or strict polynomial inequalities.}.
One can show that the projection of a semi-algebraic set (keeping a smaller number of coordinates)  is still semi-algebraic  (Tarski-Seidenberg elimination theroem, see for instance Benedetti and Risler), so  $A^*=\{(\lambda, v_\lambda), \lambda \in (0,1]\}$ is also a semi-algebraic subset of $\R\times \R^K$. This implies the existence of  a bounded Puiseux series development of $v_\lambda$ in a neighborhood of  $\lambda=0$.

\begin{thm} \label{thm1} (Bewley Kohlberg) There exists  $\lambda_0>0$, a positive integer $M$, coefficients $ r_m  \in \R^K$ for each $m\geq 0$ such that for all  $\lambda \in (0,\lambda_0]$, and all $k$ in $K$:

$$v_\lambda(k)= \sum_{m=0}^{\infty} r_m(k) \;  \lambda^{m/M}.$$

\end{thm}

So when $\lambda$ is close to 0,  for each $k$ $v_\lambda(k)$   is a power series of $\lambda^{1/M}$.

Example 1 :  $v_\lambda=\frac{1-\sqrt{\lambda}}{1- \lambda}= (1-\sqrt{\lambda})(1+\lambda +  ...+ \lambda^n+....)$

\begin{cor} $\;$

1) $v_\lambda$ converges when $\lambda$ goes to 0.

2) $v_\lambda$ has bounded variation at 0, i.e. for any  sequence $(\lambda_i)_{i\geq 1}$ of discount factors decreasing to 0, we have $\sum_{i \geq 1} \|v_{\lambda_{i+1}}-v_{\lambda_i}\| < \infty$.

3) $v_n$ also converges, and $\lim_{n \to \infty} v_n= \lim_{\lambda \to 0} v_\lambda.$
\end{cor}

\noindent{\bf Proof:} 1) is clear by the Puiseux series development.

2) also comes from this development. Fix $k$ in $K$. When $\lambda$ is small enough, $v_\lambda(k)=f_k(\lambda^{1/M})$ where $f_k$ is a power series with positive radius of convergence, hence
$$\frac{\partial v_\lambda(k)}{\partial \lambda}(\lambda)=\frac{1}{M}f'_k(\lambda^{1/M}) \lambda^{1/M-1}.$$
\noindent so that there exists a bound $C$ such that for $\lambda$ small enough, $|\frac{\partial v_\lambda(k)}{\partial \lambda}(\lambda)|\leq C  \lambda^{1/M-1}.$ Now, if $0<\lambda_2<\lambda_1$, 
$|v_{\lambda_1}(k)-v_{\lambda_2}(k)|\leq \int_{\lambda_2}^{\lambda_1} C  \lambda^{1/M-1}d\lambda= CM (\lambda_1^{1/M}-\lambda_2^{1/M}),$
and the result follows.

3) The idea\footnote{The following proof is, I believe,  due to A. Neyman.} is to compare $v_n$ with the value $w_n:=v_{1/n}$ of the $\frac{1}{n}$ discounted game.  Using the Shapley operator, we have for all $n$:
$$v_{n+1}=\frac{1}{n+1}\Psi(nv_n),\; {\rm and}\; w_{n+1}=\frac{1}{n+1}\Psi(nw_{n+1}).$$
\noindent Since $\Psi$ is non expansive,
$\|w_{n+1}-v_{n+1}\|\leq \frac{n}{n+1}\|w_{n+1}-v_n\|\leq \frac{n}{n+1}(\|w_{n+1}-w_n\| + \|w_{n}-v_n\|)$. We obtain:
$$(n+1)\|w_{n+1}-v_{n+1}\|-n \|w_n-v_n\|\leq n \|w_{n+1}-w_n\|.$$
\noindent And summing these inequalities from $n=1$ to $m$ gives:
$$\|w_{m+1}-v_{m+1}\|\leq \frac{1}{m+1}\sum_{n=1}^m n \|w_{n+1}-w_n\|.$$
It is a simple exercise to show that if $(a_n)_n$ is a sequence of non negative real numbers satisfying $\sum_{n=1}^\infty a_n<\infty$,  the sequence $(na_n)_n$ Cesaro-converges to 0.
By the bounded variation property, we have $\sum_{n=1}^\infty    \|w_{n+1}-w_n\|<\infty$. We conclude that $\|w_{m+1}-v_{m+1}\|\ \xrightarrow[m \to \infty]{} 0.$\\

 Bewley and Kohlberg also provided an example where $v_n$ is equivalent to  $\frac{\ln n}{n}$ when $n$ goes to infinity, so $(v_n)$ may not have a Puiseux series development in $n$.
 
 \subsection{Uniform value}
 
 We fix here the initial state $k_1$, and omit the dependance on the initial state for a while. We know that $\lim_n v_n=\lim_\lambda v_\lambda$ exists, so we approximately know the value of the stochastic game when $n$ is large and known to the players (and when $\lambda$ is small and known to the players). But this does not tell us if the players can play approximately well when they don't know themselves exactly how large is $n$ or how small is $\lambda$. Do the players have nearly optimal strategies that are robust with respect to the time horizon or the discount factor ? This property is captured by the notion of {\it uniform value},   which may be considered as the nectar of stochastic games.
 
 \begin{defi}  \label{def3}Let $v$ be a real.
 
 Player 1 can uniformly guarantee  $v$ in the stochastic game if: $\forall \varepsilon  >0$, $\exists \sigma\in \Sigma$, $\exists n_0$, $\forall  n \geq n_0$, $\forall \tau \in {\cal T}$, $\gamma_n(\sigma, \tau) \geq v-\varepsilon.$  

Player 2 can  uniformly guarantee $v$ in the stochastic game if: $\forall \varepsilon  >0$, $\exists \tau \in {\cal T}$, $\exists n_0$, $\forall n \geq n_0$, $\forall \sigma \in \Sigma$, $\gamma_n(\sigma, \tau) \leq v+\varepsilon.$

If $v$ can be uniformly guaranteed by both players, then $v$ is called the uniform value of the stochastic game.

 \end{defi}

It is easily shown that the uniform value, whenever it exists, is unique. The largest quantity uniformly guaranteed by Player 1, resp. smallest quantity uniformly guaranteed by Player 2, can be denoted by:
$$\underline{v}= \sup_{\sigma} \liminf_n \left( \inf_\tau \gamma_n( \sigma, \tau)\right), \; \overline{v}=\inf_\tau  \limsup_n \left( \sup_{\sigma} \gamma_n( \sigma, \tau)\right).$$
Plainly, $\underline{v}\leq \lim_n v_n \leq  \overline{v}.$ The uniform value exists if and only if $\underline{v}= \overline{v}$. Whenever it exists it is equal to $\lim_n v_n=\lim_\lambda v_\lambda$, and for each $\varepsilon>0$ there exists $\lambda_0>0$, $\sigma$ and $\tau$ such that for all $\lambda\leq \lambda_0$, $\sigma'$ and $\tau'$ we have:
$\gamma_\lambda(\sigma, \tau') \geq v-\varepsilon$ and $\gamma_\lambda(\sigma', \tau) \leq v+\varepsilon.$

\subsubsection{The Big Match}

The Big Match is the absorbing stochastic game described by:

\begin{center}
 $  
\begin{array}{cc}
 \; & \begin{array}{cc}
 L\;\; & \; \; \; R \\
 \end{array} \\
 \begin{array}{c}
 T \\
 B \\
 \end{array} &
 \left( \begin{array}{cc}
  1^* & 0^* \\
0 & 1 \\
\end{array}\right)\\
\end{array} $
\end{center}

\noindent It was introduced by Gillette in 1957. We have seen that $\lim v_n =\lim v_\lambda=1/2$ here. It is easy to see that player 2 can uniformly guarantee 1/2 by playing at each stage the mixed action 1/2 $L$ +1/2 $R$ independently of everything. It is less easy to see what can be uniformly guaranteed by player 1, and one can show that no stationary or Markov strategy of Player 1 can uniformly guarantee a positive number here. However, Blackwell and Ferguson (1968) proved   that the uniform value of the Big Match exists.

\begin{pro} \label{pro2} The Big Match has a uniform value
\end{pro}

\noindent{\bf Proof:} All we have to do is prove that Player 1 can uniformly guarantee $1/2-\varepsilon$ for each $\varepsilon>0$. First define the following random variables, for all positive integer $t$: $g_t$ is the payoff of player 1 at stage $t$, $i_t\in \{T,B\}$ is the action played by player 1 at stage $t$, $j_t\in \{L,R\}$ is the action played by player 2 at stage $t$, 
$L_t=\sum_{s=1}^{t-1}  {\bf 1}_{j_s=L}$ is the number of stages in 1,...,$t-1$ where player 2 has played $L$,  $R_t=\sum_{s=1}^{t-1}  {\bf 1}_{j_s=R}=t-1-L_t$ is the number of stages in 1,...,$t-1$ where player 2 has played $R$, and $m_t=R_t-L_t\in \{-(t-1),...,0,...,t-1\}$. $R_1=L_1=m_1=0$.

Given a fixed parameter $M$ (a positive integer) let us define  the following strategy $\sigma_M$ of player 1: at any stage $t$, $\sigma_M$ plays $T$ with probability $\frac{1}{(m_t+M+1)^2}$, and $B$ with the remaining probability.

Some intuition for $\sigma_M$ can be given. Assume we are  still in the non  absorbing state at stage $t$. If player 2 has played $R$ often at past stages, player 1 is doing well and has received good payoffs,    $m_t$ is large and $\sigma_M$ plays the risky action $T$ with small probability. On the other hand if Player 2 is playing $L$ often,  player 1 has received low payoffs but Player 2 is taking high risks;     $m_t$ is small  and $\sigma_M$ plays the risky action $T$ with high probability.  

Notice that $\sigma_M$ is well defined. If $m_t=-M$ then $\sigma_M$ plays $T$ with probability 1 at stage $t$ and then the game is over. So the event $m_t \leq -M-1$ has probability 0 as long as the play is in the non absorbing state. At any stage $t$ in the non absorbing state, we have $-M\leq m_t \leq t-1$, and $\sigma_M$ plays $T$ with a probabilty in the interval $[\frac{1}{(M+t)^2}, 1].$

We will show that $\sigma_M$ uniformly guarantees $\frac{M}{2(M+1)}$, which is close to 1/2 for $M$ large. More precisely we  will prove that:

\begin{equation} \label{eq1}\forall T\geq 1,  \forall M \geq 0, \forall \tau \in {\cal T}, \; \E_{\sigma_M,\tau}\left(\frac{1}{T} \sum_{t=1}^T g_t\right) \geq \frac{M}{2(M+1)}-\frac{M}{2T}.\end{equation}

To conclude the proof of proposition \ref{pro2}, we now prove (\ref{eq1}). Notice that  we can restrict attention to   strategies of player 2 which are pure, and (because there is a unique relevant  history of moves of player 1) independent of the history. That is, we can assume w.l.o.g. that player 2 plays a fixed deterministic  sequence $y=(j_1,...j_t,...)\in \{L,R\}^\infty.$  

$T$ being fixed until the end of the proof, we define the random variable $t^*$ as the time of absorption: $$t^*=\inf\{s\in\{1,...,T\}, i_s=T\},\; {\rm with \; the \; convention\;}t^*=T+1 \; {\rm if}\; \forall s\in\{1,...,T\}, i_s=B$$

Recall that $R_t=m_t+L_t=t-1-L_t$, so that $R_t=\frac{1}{2} (m_t+t-1)$. For $t\leq t^*$, we have $m_t\geq -M$, so:
$$R_{t^*}\geq \frac{1}{2} (t^*-M-1)$$

Define also $X_t$ as the following fictitious payoff of player 1: $X_t=1/2$ if $t\leq t^*-1$, $X_t=1$ if $t\geq t^*$ and $j_{t^*}=L$, and $X_t=0$ if $t\geq t^*$ and $j_{t^*}=R$. $X_t$ is the random variable of the limit value of the current state.

A simple computation shows:

\begin{eqnarray*}
\E_{\sigma_M,y} \left(\frac{1}{T} \sum_{t=1}^T g_t\right) & = & \E_{\sigma_M,y} \; \frac{1}{T}  (R_{t^*}+(T-t^*+1) {\bf 1}_{j_{t^*}=L})\\
 & \geq &  \E_{\sigma_M,y} \; \frac{1}{T}  (\frac{1}{2} (t^*-M-1)+(T-t^*+1) {\bf 1}_{j_{t^*}=L})\\
 & \geq  & -\frac{M}{2T} +  \E_{\sigma_M,y} \; \frac{1}{T}  (\frac{1}{2} (t^*-1)+(T-t^*+1) {\bf 1}_{j_{t^*}=L})\\
 & \geq & -\frac{M}{2T} +  \E_{\sigma_M,y} \; \left( \frac{1}{T} \sum_{t=1}^T X_t\right)
\end{eqnarray*}

To prove (\ref{eq1}), it is thus enough to show the following lemma.

\begin{lem} \label{lem1} For all $t$ in $\{1,...,T\}$, $y$ in $\{L,R\}^\infty$ and $M$ in $\N$, $ \E_{\sigma_M,y} \;  (X_t) \geq \frac{M}{2(M+1)}$.
\end{lem}

\noindent Proof of the lemma. The proof is by induction on $t$. For $t=1$, 
$ \E_{\sigma_M,y} \;  (X_1)=\frac{1}{2} ( 1- \frac{1}{(M+1)^2})+ \frac{1}{(M+1)^2} \1_{j_1=L} $ $\geq$ $\frac{1}{2} ( 1- \frac{1}{(M+1)^2})$ $\geq $ $ \frac{M}{2(M+1)}$.
  
  Assume the lemma is true for $t\in \{1,...,T-1\}$. Consider $y=(j_1,...)$ in $\{L,R\}^\infty$, and write $y=(j_1,y_+)$ with $y_+=(j_2,j_3,...)\in \{L,R\}^\infty$.
  If $j_1=L$, $ \E_{\sigma_M,y} \;  (X_{t+1})=\frac{1}{(M+1)^2}1 +  ( 1- \frac{1}{(M+1)^2}) \E_{\sigma_{M-1},y_+}(X_t)$. By the induction hypothesis, $ \E_{\sigma_{M-1},y_+}(X_t)\geq  \frac{M-1}{2M}$, so $ \E_{\sigma_M,y} \;  (X_{t+1})\geq  \frac{M}{2(M+1)}$. Otherwise $j_1=R$, and $ \E_{\sigma_M,y} \;  (X_{t+1})=   ( 1- \frac{1}{(M+1)^2}) \E_{\sigma_{M+1},y_+}(X_t)$ $\geq  ( 1- \frac{1}{(M+1)^2})\frac{M+1}{2(M+2)}$ $= \frac{M}{2(M+1)}$. The lemma is proved, concluding the proof of proposition \ref{pro2}.

\subsubsection{The existence result}

The following theorem is due to J-F. Mertens and A. Neyman (1982).

\begin{thm}\label{thm2} (Mertens Neyman)

Every zero-sum stochastic game  with finitely many states and actions has a uniform value.
\end{thm}

The rest of this section is devoted to the proof of theorem \ref{thm2}. Without loss of generality we assume that all payoffs are in $[0,1]$, and fix $\varepsilon \in (0,1)$ in the sequel.

We know by the algebraic approach that there exists $C>0$,   $M\geq 1$, $\lambda_0>0$ such that for all $0<\lambda_1<\lambda_2\leq \lambda_0$:
$$\|v_{\lambda_1}-v_{\lambda_2}\|\leq \int_{\lambda_1}^{\lambda_2} \psi(s)ds\;\;{\rm with\;}\; \psi(s)=\frac{C}{s^{1-1/M}}.$$
\noindent All is needed about  $\psi$ is that it is non negative  and integrable: $ \int_{0}^{1} \psi(s)ds<\infty.$  

\begin{defi}\label{def4} Define the mapping $D$ from $(0, \lambda_0]$ to $\R$ by:
$$D(y)=\frac{12}{\varepsilon} \int_{y}^{\lambda_0} \frac{\psi(s)}{s} ds + \frac{1}{\sqrt{y}}.$$
\end{defi}
The proof of the next lemma is left to the reader.
\begin{lem}\label{lem2} $\;$

a) $D$ is positive,   decreasing, $D(y) \xrightarrow[y \to0]{} +\infty$ and $\int_{0}^{\lambda_0} D(y)dy<\infty.$

b) $D(y(1-\varepsilon/6))-D(y) \xrightarrow[y \to0] {}+\infty$ and $D(y)-D(y(1+\varepsilon/6))  \xrightarrow[y \to0] {}+\infty$.

c) $yD(y) \xrightarrow[y \to0] {}0$.
\end{lem}

\begin{defi}\label{def5} Define the mapping $\varphi$ from $[0, \lambda_0]$ to $\R$   by:
$$\varphi(\lambda)=\int_{0}^\lambda D(y)dy- \lambda D(\lambda).$$\end{defi}


Note that $\varphi$ is increasing and $\varphi(0)=\lim_{\lambda \to 0} \varphi(\lambda)=0$. \\

We  fix the initial state $k_1$ and denote the limit value $\lim_\lambda v_\lambda(k_1)$ by $v(k_1)$. We now define a nice  strategy $\sigma $ for player 1 in the stochastic game with initial state $k_1$.  While playing at some stage $t+1$, player 1 knows the current state $k_{t+1}$ and the previous payoff $g_t$, he will update a fictitious discount factor $\lambda_{t+1}$ and play at stage $t+1$ a stationary optimal strategy in the stochastic game with discount factor $\lambda_{t+1}$ and initial state $k_{t+1}$. The definition of the sequence of random discount factors $(\lambda_t)_t$ below, joint with  the introduction of an auxiliary sequence  $(d_t)_t$, will end the definition of $\sigma$.

One first chooses $\lambda_1>0$ such that: (i) $v_{\lambda_1}(k_1)\geq v(k_1)-\varepsilon$, (ii) $\varphi(\lambda_1)<\varepsilon$, and (iii) $\forall y \in (0, \lambda_1],\;D(y(1-\varepsilon/6))-D(y)>6$ and  $D(y)-D(y(1+\varepsilon/6))>6$. Put $d_1=D(\lambda_1)$, and by induction define for each $t\geq 1$:
$$d_{t+1}=\max\{d_1,d_t+g_t-v_{\lambda_t}(k_{t+1})+4 \varepsilon\},\; {\rm and}\;\lambda_{t+1}=D^{-1}(d_{t+1}).$$

We have  $\lambda_{t+1}\leq \lambda_1$ for each $t$. Notice  that if the current payoff $g_t$ is high, then $\lambda_{t+1}$ will have a tendency to be small : player 1 plays in a patient way. On the contrary if $g_t$ is small then $\lambda_{t+1}$ will have a tendency to be large  : player 1 plays   for the short-run payoffs.
$\sigma$ being defined, we now fix a strategy $\tau$ of player 2. We simply write $\P$ for $\P_{k_1, \sigma, \tau}$ and $\E$ for $\E_{k_1,\sigma, \tau}$.\\

By construction, the following properties hold on every   play. The proofs of a), b) and d) are left to the reader.
\begin{lem}\label{lem3} $\;$ For all $t\geq 1$,

a) $|d_{t+1}-d_t|\leq 6$,

b)  $|\lambda_{t+1}-\lambda_t|\leq \frac{\varepsilon \lambda_t}{6}$,

c) $|v_{\lambda_t}(k_{t+1})-v_{\lambda_{t+1}}(k_{t+1})|\leq \varepsilon \lambda_t$.

d) $d_{t+1}-d_t\leq g_t-v_{\lambda_t}(k_{t+1})+ 4 \varepsilon + \1_{\lambda_{t+1}=\lambda_1}.$

\end{lem}

\noindent{\bf Proof of } $c)$:  \begin{eqnarray*}
|v_{\lambda_t}(k_{t+1})-v_{\lambda_{t+1}}(k_{t+1})| & \leq & \|v_{\lambda_t}-v_{\lambda_{t+1}}\|\\
 & \leq & \left| \int_{\lambda_t}^{\lambda_{t+1}} \psi(s) ds\right|\\
  & \leq & \max\{\lambda_t, \lambda_{t+1}\}  \left| \int_{\lambda_t}^{\lambda_{t+1}} \frac{\psi(s)}{s} ds\right|\\
    & \leq & 2 \lambda_t  \left| \int_{\lambda_t}^{\lambda_{t+1}} \frac{\psi(s)}{s} ds\right|.
 \end{eqnarray*}
 Now, $$ \int_{\lambda_t}^{\lambda_{t+1}} \frac{\psi(s)}{s}ds= \frac{\varepsilon}{12} \left(D(\lambda_t)- \frac{1}{\sqrt{\lambda_t}}\right) -  \frac{\varepsilon}{12} \left(D(\lambda_{t+1})- \frac{1}{\sqrt{\lambda_{t+1}}}\right)=  \frac{\varepsilon}{12}\left((d_t-d_{t+1})+ (\frac{1}{\sqrt{\lambda_{t+1}}}-\frac{1}{\sqrt{\lambda_t}})\right).$$
 If $\lambda_t \leq \lambda_{t+1}$, $0\leq  \int_{\lambda_t}^{\lambda_{t+1}} \frac{\psi(s)}{s}ds\leq   \frac{\varepsilon}{2}$ by point a) of lemma \ref{lem3}. So  $  \left| \int_{\lambda_t}^{\lambda_{t+1}} \frac{\psi(s)}{s} ds\right|\leq  \frac{\varepsilon}{2}$, and this also holds if $\lambda_t \geq \lambda_{t+1}$. We obtain $|v_{\lambda_t}(k_{t+1})-v_{\lambda_{t+1}}(k_{t+1})|\leq \varepsilon \lambda_t$.\\

\begin{defi}\label{def6} Define the random variable $$Z_t=v_{\lambda_t}(k_t)-\varphi(\lambda_t).$$
\end{defi}
When $\lambda_t$ is close to 0, $Z_t$ is close to $v(k_t)$.

\begin{pro}\label{pro3} $\;$ $(Z_t)_t$ is a sub-martingale, and for all $t\geq 1$:
$$\E(Z_{t})\geq 2 \varepsilon \E(\sum_{s=1}^{t-1} \lambda_s)+ Z_1.$$
\end{pro}

Proposition \ref{pro3} is the key to Mertens and Neyman's proof. Assume for the moment the proposition  and let us see how the proof of the theorem follows. \\

We have for each $t\geq 1$, $\E(Z_t)\geq Z_1$, so $\E(v_{\lambda_t}(k_t))\geq v_{\lambda_1}(k_1)-\varphi(\lambda_1)+ \E(\varphi(\lambda_t))\geq v_{\lambda_1}(k_1)-\varphi(\lambda_1)$, so 

\begin{equation} \label{eq2} \E(v_{\lambda_t}(k_t))\geq v_{\lambda_1}(k_1)-\varepsilon.\end{equation}
Since $Z_{t+1}\leq 1$, we have by proposition  \ref{pro3} that $2 \varepsilon \E(\sum_{s=1}^{t} \lambda_s)\leq 1-Z_1\leq 1+ \varepsilon$, so $ \E(\sum_{s=1}^{t} \lambda_s)\leq \frac{1}{\varepsilon}.$ We obtain  $ \E(\sum_{s=1}^{t}\lambda_1 \1_{\lambda_1= \lambda_s})\leq \frac{1}{\varepsilon}$, and 
\begin{equation} \label{eq3} \E (\sum_{s=1}^{t} \1_{\lambda_1= \lambda_s})\leq \frac{1}{\lambda_1\varepsilon}.\end{equation}

Using d) and c)  of lemma \ref{lem3}, we have: 
$$g_t\geq    v_{\lambda_{t+1}}(k_{t+1}) - \varepsilon \lambda_t + (d_{t+1}-d_t)- 4 \varepsilon - \1_{\lambda_{t+1}=\lambda_1}.$$
\noindent so for each $T$, $$\E\left(\frac{1}{T}\sum_{t=1}^T g_t\right) \geq
   \E\left(\frac{1}{T}v_{\lambda_{t+1}}(k_{t+1})\right) - \varepsilon \E\left(\frac{1}{T}\sum_{t=1}^T \lambda_t\right) + \E\left(\frac{1}{T}(d_{T+1}-d_1)\right) - 4 \varepsilon - \frac{1}{T} \E\left(\sum_{t=1}^T  \1_{\lambda_{t+1}=\lambda_1}\right)$$
Unsing the inequalities (\ref{eq2}) and (\ref{eq3}), we obtain    $$\E\left(\frac{1}{T}\sum_{t=1}^T g_t\right)  \geq v_{\lambda_1}(k_1)-\varepsilon -\varepsilon -\frac{d_1}{T}- 4 \varepsilon -\frac{1}{\varepsilon \lambda_1 T}.$$
 And for $T$ large enough, we have:
  $$\E\left(\frac{1}{T}\sum_{t=1}^T g_t\right) \geq v(k_1) -8 \varepsilon,$$
  independently of the strategy $\tau$ of player 2. This shows that player 1 uniformly guarantees $v(k_1)$ in the stochastic game with initial state $k_1$. By symmetry, player 2 can do as well and theorem \ref{thm2} is proved.\\
  
  We finally come back to the proof of the key proposition.\\

  \noindent{\bf Proof of proposition \ref{pro3}}: Fix $t\geq 1$, and define
 $C_1= \varphi(\lambda_t)-\varphi(\lambda_{t+1})$, $C_2=v_{\lambda_{t+1}}(k_{t+1})-v_{\lambda_t}(k_{t+1})$ and $C_3=\lambda_t(g_t-v_{\lambda_t}(k_{t+1}))$. A simple computation shows that:
 $$Z_{t+1}-Z_t- (C_1+C_2-C_3)= \lambda_t g_t+(1-\lambda_t) v_{\lambda_t}(k_{t+1})-v_{\lambda_t}(k_t).$$
Denote by ${\cal H}_t$ the $\sigma$-algebra generated by histories in $(K \times I \times J)^{t-1}\times K$ (before players play at stage $t$), by definition of $\sigma$ one has:
$$\E(\lambda_t g_t +(1-\lambda_t) v_{\lambda_t}(k_{t+1})| {\cal H}_t)\geq v_{\lambda_t}(k_t).$$
   Consequently, we obtain:
\begin{equation} \label{eq4} \E(Z_{t+1}-Z_t| {\cal H}_t)\geq \E(C_1+C_2-C_3 |  {\cal H}_t).\end{equation}
   We have $|C_2|\leq \varepsilon \lambda_t$ by point c) of lemma \ref{lem3}.
 By definition of $d_{t+1}$, we have $d_{t+1}-d_t\geq g_t-v_{\lambda_t}(k_{t+1})+ 4 \varepsilon$, hence $C_3\leq \lambda_t (d_{t+1}-d_t)-4 \varepsilon \lambda_t$. We now prove
   
   \begin{equation}\label{eq5} C_1\geq \lambda_t(d_{t+1}-d_t)- \varepsilon \lambda_t\end{equation}
 
 \noindent If $\lambda_{t+1}<\lambda_t$, then $d_{t+1}>d_t$ and $C_1= \varphi(\lambda_t)-\varphi(\lambda_{t+1})$ $\geq$ $ \lambda_t( d_{t+1}-d_t) -(\lambda_t-\lambda_{t+1})(d_{t+1}-d_t)$ $\geq$  $ \lambda_t( d_{t+1}-d_t) - \varepsilon \lambda_t$ by a) and b) of lemma \ref{lem3}. If $\lambda_{t+1}>\lambda_t$, then  $d_{t+1}<d_t$ and $\varphi(\lambda_{t+1})-\varphi(\lambda_{t})\leq \lambda_{t+1} (d_t-d_{t+1})$ $=$ $ \lambda_t  (d_t-d_{t+1})+ (\lambda_{t+1}-\lambda_t)  (d_t-d_{t+1})$ $\leq$ $ \lambda_t(  d_{t}-d_{t+1})+ \varepsilon \lambda_t$. And (\ref{eq5}) is proved.\\
   
 Back to inequality (\ref{eq4}), we obtain: $$\E(Z_{t+1}-Z_t| {\cal H}_t)\geq \E\left(\lambda_t(d_{t+1}-d_t)- \varepsilon \lambda_t - \varepsilon \lambda_t - \lambda_t (d_{t+1}-d_t)+4 \varepsilon \lambda_t  | {\cal H}_t\right)=2 \varepsilon  \E\left(\lambda_t| {\cal H}_t)\right.$$
 \noindent which proves that $(Z_t)_t$ is a sub-martingale  and for all $t\geq 0$, $\E(Z_{t+1})\geq 2 \varepsilon \E(\sum_{s=1}^{t} \lambda_s)+ Z_1.$ This ends the proof of proposition \ref{pro3}.
 
  \newpage
 
 \section{A few  extensions and recent results}
 
 We want to go beyond the ``simple" case of finitely many states and actions. It is interesting to start with the one-player case, which is fairly understood.
 
 \subsection{1-Player games}
 
 \subsubsection{General results: the long-term value}
 
 We consider a general dynamic programming problem with bounded payoffs:
 \textcolor{black}{  $\Gamma(z_0)=(Z,F,r,z_0)$} given by  a non empty set of states $Z$, an  initial state $z_0$,  a transition correspondence $F$ from $Z$ to $Z$ with non empty values, and a reward mapping $r$ from $Z$ to $[0,1]$.   Here $Z$ can be any set, and for each state $z$ in $Z$, $F(z)$ is a non empty subset of $Z$. (An equivalent MDP variant of the model exists with an explicit set of actions $A$, and transitions given by a function from $Z\times A$ to $Z$.)
 
   \vspace{0,2cm}
   
A player  chooses $z_1$ in $F(z_0)$, has a payoff of   $r(z_1)$, then he chooses $z_2$ in $F(z_1)$, etc...  

The set of admissible plays at $z_0$ is defined  as: $S(z_0)=\{s=(z_1,...,z_t,...)\in Z^{\infty}, \forall t \geq 1, z_t \in F(z_{t-1})\}.$
  
  \vspace{0,2cm}

For $n\geq 1$, the value of the $n$-stage problem with initial state $z$ is defined as: 
$$v_n(z)= \sup_{s \in S(z)} \gamma_n(s), \; {\rm where} \; \gamma_n(s)= \frac{1}{n} \sum_{t=1}^n r(z_t).$$

   \vspace{0,2cm}
  For $\lambda \in (0,1]$, the value of the  $\lambda$-discounted problem with initial state $z$ is defined as:
 $$v_\lambda(z)=\sup_{s \in S(z)} \gamma_\lambda(s), \; {\rm where} \; \gamma_\lambda(s)=\lambda \sum_{t=1}^{\infty} {(1-\lambda)}^{t-1} r(z_t).$$

 More generally, define  an   \textcolor{black}{evaluation       $\theta=(\theta_t)_{t \geq 1}$}  as a probability   on  positive integers. The $\theta$-payoff of a play $s=(z_t)_{t\geq 1}$ is  
$\gamma_\theta(s)=   \sum_{t=1}^\infty \theta_tr(z_t),$
  and   \textcolor{black}{the $\theta$-value of   $\Gamma(z)$ is 
$$v_\theta(z)= \sup_{s \in S(z)} \gamma_\theta(s).$$} The set of all evaluations is denoted by $\Theta$. The total variation of an evaluation $\theta$ is defined as: $TV(\theta)=\sum_{t=1}^\infty |\theta_{t+1}-\theta_t|.$
Given an evaluation $\theta=\sum_{t\geq 1} \theta_t \delta_t$ (here $\delta_t$ is the Dirac measure on stage $t$) and some $m\geq 0$,   we write ${{v}}_{m, \theta}$ for the value function associated to the shifted evaluation $\theta \oplus m = \sum_{t=1}^\infty  \theta_t \delta_{m+t}.$  \\

 What can be said in general about the convergence of $(v_n)_n$, when $n\to \infty$, of  $(v_\lambda)_\lambda$, when $\lambda \to 0$, or more generally of   $(v_{\theta^k})_k$, when $(\theta^k)_k$ is a sequence of evaluations such that $TV(\theta^k)\to_{k \to \infty} 0$? Many things, if we focus on {\it uniform  convergence.} \\
 
We now only consider uniform convergence of the value functions.  Denote by ${\cal V}$ the set of functions from $Z$ to $[0,1]$, endowed with the supremum metric $d_{\infty}(v,v')=\sup_{z \in Z} |v(z)-v(z')|$. Saying  that a sequence $(v^k)_{k \geq 1}$ of functions from $Z$ to $[0,1]$ uniformly converges  is trivially equivalent to saying that the sequence $(v^k)$ converges in the metric space ${\cal V}$. Notice  that in a metric space, convergence of a sequence $(v^k)$ happens if and only if:
 
 1) the sequence $(v^k)_k$ has at most one limit\footnote{A limit point of $(v^k)_k$ being defined as  a limit of a converging subsequence of   $(v^k)_k$.}, and 
 
 2) the set $\{v^k, k\geq 1\}$ is totally bounded\footnote{For each $\varepsilon>0$, the set can be covered by finitely many balls of radius $\varepsilon$. Equivalently, the completion of the set is compact. Equivalently, from any sequence in the set one can extract a Cauchy subsequence.}.

 The above equivalence holds for any sequence in a metric space. But here we consider the special case of value functions of a dynamic programming problem, with long term limits. It will turn out that 1) is automatically satisfied.

\begin{defi} Define for all $z$ in $Z$, \textcolor{black}{$$v^*(z)=\inf_{\theta \in \Theta} \; \sup_{m \geq 0} \;  {{v}}_{m, \theta}(z).$$} 
\end{defi}

The following results apply in particular to the sequences $(v_n)_n$ and $(v_\lambda)_\lambda$ (see Renault, 2011).
\begin{thm} (Renault, 2014)

 Consider a sequence of evaluations  $(\theta^k)_k$   such that $TV(\theta^k)\to_{k \to \infty} 0$. 

Any limit point of $(v_{\theta^k})_k$ is $v^*$.
\end{thm}

The proof is omitted.

\begin{cor} \label{cor1} Consider a sequence of evaluations  $(\theta^k)_k$   such that $TV(\theta^k)\to_{k \to \infty} 0$. 

1) If  $(v_{\theta^k})_k$  \; {\rm  converges}, the limit is $v^*$.

2) \begin{eqnarray*}
 (v_{\theta^k})_k  \; {\rm  converges} & \Longleftrightarrow & {\rm  the \; set  } \; \{v_{\theta^k}, k\geq 1\} \; {\rm is \; totally \; bounded,}\\
\; & \Longleftrightarrow &  {\rm  the \; set \; } \; \{v_{\theta^k}, k\geq 1\} \cup \{v^*\}  \; {\rm is \; compact.}
\end{eqnarray*}

3) Assume that $Z$ is endowed with a distance $d$ such that:  a) $(Z,d)$ is a  totally bounded  metric space, and b)  the family $(v_\theta)_{\theta \in \Theta}$  is   uniformly equicontinuous. Then there is general uniform convergence of the value functions to $v^*$, i.e. $$\forall \varepsilon >0, \exists \alpha >0, \forall \theta \in \Theta \; s.t. \; TV(\theta)\leq \alpha, \; \|v_\theta-v^*\| \leq \varepsilon.$$ 

4) Assume that  $Z$ is endowed with a distance $d$ such that: a)  $(Z,d)$ is a    precompact metric space, b)   $r$ is   uniformly continuous,  and c) $F$   is   non expansive, i.e. $\forall z\in Z, \forall z' \in Z, \forall z_1 \in F(z), \exists z'_1 \in F(z') \; s.t. \; d(z_1,z'_1)\leq d(z,z').$ Same conclusions as corollary 3. 
\end{cor}
  
The above results can be extended to the case of stochastic dynamic programming, (i.e. when  $F(z)$ is a set of probability distributions on $Z$ for each $z$). In this case it is often convenient to define the value functions $v_n$, $v_\lambda$, $v_\theta$ directly by their Bellman equations. 

 Notice that life is much simpler in  the particular case where the problem is leavable, i.e. when $z\in F(z)$ for each $z$. Then without any assumption, $(v_n)$ and $(v_\lambda)$ {\it pointwise} converge to $v^*$, where:
 $v^*=\inf\{v: Z\to [0,1],  excessive\footnote{$v$ excessive means that $v(z)\geq v(z')$ if $z'\in F(z)$, i.e. that $v$ is non increasing on any trajectory.}, v\geq r\}$.

  \subsubsection{The uniform convergence  of $(v_n)_n$ and $(v_\lambda)_\lambda$ are equivalent.}
    
 The results of the previous subsection show in particular  that if $(v_n)$ and $(v_\lambda)$ uniformly converge, they have the same limit.  For these two particular sequences of evaluations, we have a stronger result.
 
  \begin{thm} (Lehrer-Sorin 1992)  In a 1-player game,  $(v_n)$ converges uniformly if and only if $(v_\lambda)$ converges uniformly. In case of convergence, the limit is the same.\end{thm}

 \subsubsection{The compact non expansive case and  the uniform value}

 We have stronger results if the state space is assumed to be compact, payoffs are continuous and transitions are non expansive. We consider here a stochastic dynamic programing problem (also called Gambling House) $\Gamma=(X,F,r,x_0)$ given by:  

$\bullet$ a non empty set of states $X$, an  initial state $x_0$, 

$\bullet$  a   transition multifunction $F$ from $X$ to $Z:=\Delta_f(X)$ with non empty values, 
 
$\bullet$  and a reward mapping $r$ from $X$ to $[0,1]$.  

\noindent Here $\Delta_f(X)$ is the set of probabilities with finite support over $X$. We assume that transitions have finite support for simplicity, however  many results concerning the limit value and is characterization can   go through without this assumption. When we will study the uniform value, this assumption will be  useful to define strategies avoiding  measurability issues.

 Here a player chooses $u_1$ in $F(x_0)$, then $x_1$ is selected according to $u_1$ and yields the payoff $r(x_1)$, then the player chooses $u_2$ in $F(x_1)$,   etc...
We define as usual the $n$- stage\; value \;function:  $v_n(x_0)= \sup_{\sigma \in S(x_0) } \E_{\sigma}  \left(\frac{1}{n} \sum_{t=1}^n  r(x_t)\right),$
 where $S(x_0)=\{\sigma=(u_1,...,u_t,...)\in Z^\infty, u_1\in F(x_0), \forall t\geq 1, u_{t+1}\in F(u_t)\} .$
We define similarly the $\lambda$-discounted value $v_\lambda(z_0)$, and more generally for any evaluation $\theta$ we have the $\theta$-value $v_\theta(z_0)$.

 We assume here that $X$ is a compact metric space with metric denoted by $d$. 
 The set   $\Delta(X)$  of Borel probability measures  over $X$ is also a compact metric space (for the weak-* topology), and we will use 
  the Kantorovich-Rubinstein metric\footnote{In the second expression, $\Pi(u,u')$ denotes the set of probabilities on $X\times X$ with first marginal $u$ and second marginal $u'$.}: for $z$ and $z'$ in $\Delta(X)$, 
 \begin{eqnarray*}d_{KR}(u,u')& =& \sup_{f: X \to \R,  1-Lip} \left|\int_{x \in X }  f(x)du(x) - \int_{x \in X }  f(x)du'(x)\right|\\
 & =&\min_{\pi\in \Pi(u,u')}\int_{(x,x')\in X \times X} d(x,x') d\pi(x,x').
 \end{eqnarray*}
 
$X$ is now viewed as a subset of $\Delta(X)$, and we assimilate an element $x$ in $X$ with the corresponding Dirac measure in $\Delta(X)$. The Graph of $\Gamma$ can be viewed as a subset of $\Delta(X)\times \Delta(X)$, and we denote by $\overline{\conv}{\rm Graph}(\Gamma)$ its closed convex hull in $\Delta(X)\times \Delta(X)$. We define the set of {\it invariant  measures} as:
$$R=\{u \in \Delta(X), (u,u)\in \overline{\conv}{\rm Graph}(\Gamma)\}$$

We will assume  that $r$ is continuous, and extend $r$ to a continuous  affine function defined on $\Delta(X)$: for $u$ in $\Delta(X)$, $r(u)$ is the expectation of $r$ with respect to $u$.
 We will also  assume non expansive transitions.
 $$\forall x\in X, \forall x'\in X, \forall u \in \Gamma(x), \exists u'\in \Gamma(x'),\; s.t. \;  d_{KR}(u,u')\leq d(x,x').$$
This assumption is always satisfied when $X$ is finite\footnote{if $d(x,x')=2$ for all $x$, $x'$, then $d_{K,R}(z,z')=\|z-z'\|_1$ for all $z$, $z'$ in $\Delta(X)$.}, or when $X$ is a simplex  and $\Gamma(x)$ is the set of splittings at $x$, i.e. the set of Borel probabilities  on $X$ with mean $x$. \\

  One can apply here a variant of property 4) of corollary \ref{cor1} to prove uniform   convergence of $(v_n)$ and $(v_\lambda)$, but we can obtain a stronger result with a better characterization of the limit value  and   the existence of the uniform value.

   \begin{thm}(Renault-Venel 2017) Assume the state space is compact, payoffs are continuous and transitions are non expansive.
   Then $(v_n)$ and $(v_\lambda)$ uniformly converge to $v^*$, where for each initial state $x$,
    $$   v^*(x)=\inf \large\{ 
w(x),  w:\Delta(X)\rightarrow [0,1] \;  {\rm affine}\;  C^0 \; s.t. $$
$$ \hspace{2cm} \textcolor{black}{(1)\;  \forall x'\in X, w(x')\geq \sup_{u \in F(x')} w(u)}  $$
$$ \hspace{1cm}  \textcolor{black}{ (2) \;  \forall u \in R, w(u)\geq r(u)   \}.}$$
Moreover,  the uniform value exists if $F$ has convex  values (or if one allows the player to play a behavior strategy, i.e. to select randomly an element $u$ in $F(x)$ while at state $x$).   \end{thm}

 The theorem also extend to general sequences of evaluations with vanishing total variation.
 
For partially observable Markov decision processes (POMDP)  with finite set of states, actions and signals, the existence of the uniform  value was first proved by Rosenberg, Solan and Vieille (2002). The present theorem can not be applied as is in this case, because transitions are not non expansive with respect to the $KR$-metric. However, an alternative metric introduced in (Renault  Venel 2017) can be used to apply the theorem to this class of games.
 
Recently, Venel and  Ziliotto (2016) proved for these models the existence of the uniform value in pure strategies, i.e. without the assumption that $F$ has convex values.

\subsection{A simple stochastic game with compact action sets and  no limit value}

Let us consider stochastic games with finitely many states, compact actions sets and continuous payoffs and transitions.   Bolte, Gaubert and Vigeral (2015) proved the existence of the limit value under some semi-algebraic (or more generally definability) conditions, and Vigeral (2013) showed that the limit value may not exist without  semi-algebraic conditions. The  simple example given here is a variant of an example from  Ziliotto, 2016 (a variant is also mentioned in Sorin Vigeral  2015). Contrary to Vigeral's example, the transitions are very simple (polynomial)  but the action set of player 1 will be a countably infinite  compact subset of the real line. 

\vspace{0,5cm}

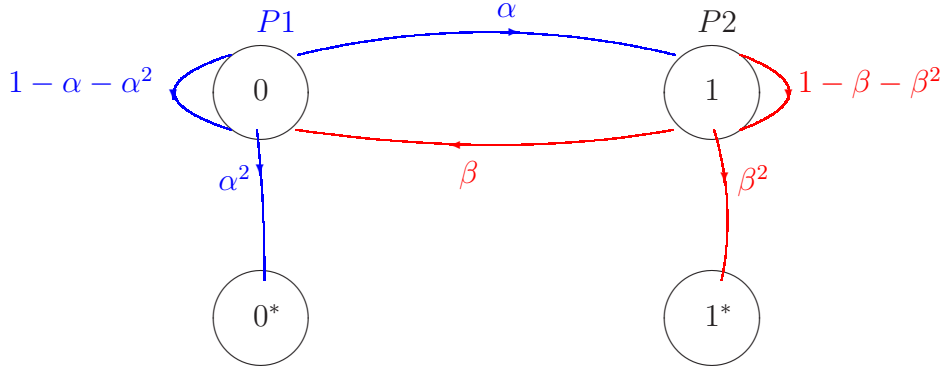
 \begin{center}
\setlength{\unitlength}{1mm}
\begin{picture}(100,60)

 \put(20,20){\circle{12}}
 \put(80,20){\circle{12}}
  \put(80,50){\circle{12}}
   \put(20,50){\circle{12}}
   
    \put(19,19){$0^*$}
     \put(79,19){$1^*$} 
       \put(19,49){$0$}
     \put(79,49){$1$} 
     
       \put(18,58){   \textcolor{blue} {  $P1$}}
          \put(78,58){$P2$}
          
             \put(50,60){   \textcolor{blue} {  $\alpha$}}
                       \put(45,38){   \textcolor{red} {  $\beta$}}

               \put(-15,50){   \textcolor{blue} {  $1-\alpha-\alpha^2$}}
                       \put(90,50){   \textcolor{red} {  $1-\beta-\beta^2$}}

               \put(13,37){   \textcolor{blue} {  $\alpha^2$}}
                       \put(82,37){   \textcolor{red} {  $\beta^2$}}
          
    \textcolor{blue} {      { \qbezier(25,55)(50,61)(75,55)}}
 \textcolor{blue} {   { \put(52,58){\vector(1,0){1}}}}

         \textcolor{red} {      { \qbezier(72,45)(50,41)(22,45)}}
 \textcolor{red} {   { \put(42,43){\vector(-1,0){1}}}}   
 
     \textcolor{blue} {      { \qbezier(14,45)(15,35)(15,25)}}
 \textcolor{blue} {   { \put(13,40){\vector(0,-1){1}}}}

      \textcolor{red} {      { \qbezier(72,45)(75,35)(73,25)}}
 \textcolor{red} {   { \put(72,39){\vector(0,-1){1}}}}

     \textcolor{blue} {      { \qbezier(5,55)(-10,50)(5,45)}}
 \textcolor{blue} {   { \put(-4,50){\vector(0,-1){1}}}}

      \textcolor{red} {      { \qbezier(70,55)(83,50)(70,45)}}
 \textcolor{red} {   { \put(75,50){\vector(0,-1){1}}}} 
 
\end{picture}
    
      \end{center}
 
 There are 4 states: $K=\{k_0,k_1,0^*,1^*\}$. States $0^*$ and $1^*$ are absorbing, and the payoff in state $k_0$, resp. $k_1$, is 0, resp. 1. In state $k_0$, Player 1 chooses $\alpha$ in some fixed set $I\subset [0,1/2]$, and the next state is $k_1$ with probability $\alpha$, $0^*$ with probability $\alpha^2$ and $k_0$ with the remaining probability $1-\alpha-\alpha^2$. Similarly, in state $k_1$ player 2 chooses $\beta$ in   $J$, and the next state is $k_0$ with probability $\beta$, $1^*$ with probability $\beta^2$ and $k_1$ with the remaining probability. To obtain divergence of the values, we introduce  a dissymmetry between players and  assume that: $$I=\{\frac{1}{2^{2n}}, n \geq 1\}, \; {\rm and }\; J=[0,1/2].$$
 
 Given a discount factor $\lambda$ in $(0,1]$, we denote by $x_\lambda$ the value of the stochastic game with discount $\lambda$ and initial state $k_0$, and by  $y_\lambda$ the value of the stochastic game with discount $\lambda$ and initial state $k_1$. The Shapley equations read:
  $$x_\lambda=\max_{\alpha \in I} (1-\lambda) ((1-\alpha-\alpha^2)x_\lambda + \alpha y_\lambda),\; {\rm and }\; y_\lambda=\min_{\beta \in J} \left(\lambda+ (1-\lambda) ((1-\beta-\beta^2)y_\lambda + \beta x_\lambda + \beta^2)\right).$$
 Substracting by $x_\lambda$ in the first equation and by $y_\lambda$ in the second equation we obtain:
 \begin{eqnarray}
 \lambda x_\lambda &= & (1-\lambda) \max_{\alpha \in I} \left(\alpha (y_\lambda-x_\lambda) - \alpha^2 x_\lambda\right) \label{eq6} \\
 \lambda y_\lambda &= &\lambda + (1-\lambda) \min_{\beta \in J} (\beta(x_\lambda-y_\lambda)+ \beta^2(1-y_\lambda)) \label{eq7}
 \end{eqnarray}
 
 Since $x_\lambda> 0$, eq. (\ref{eq6}) gives that $y_\lambda>x_\lambda$.

 \begin{lem}\label{lem4} For $\lambda\leq 1/5$,
\begin{equation}\label{eq8}4 \lambda {(1-y_\lambda)}^2=(1-\lambda){(y_\lambda-x_\lambda)}^2.\end{equation}
 \end{lem}
 \noindent{\bf Proof:} Consider the convex minimization problem of player 2 given by (\ref{eq7}). Either the minimum is achieved for the interior point $\beta_\lambda:=\frac{y_\lambda-x_\lambda}{2(1-y_\lambda)}$, and we get equation (\ref{eq8}). Or in the case when $y_\lambda-x_\lambda>1-y_\lambda$, it is achieved for $\beta=1/2$ and we obtain $(1-y_\lambda)(1+3 \lambda)=2 (1-\lambda)(y_\lambda-x_\lambda).$ So in this case $1+3\lambda>2(1-\lambda)$, which is not possible for $\lambda \leq 1/5$.\\
 
Equation (\ref{eq8}) will be our main tool.  An immediate consequence   is that $y_\lambda-x_\lambda\longrightarrow_{\lambda \to 0}0$. Look now at the concave maximization problem of player 1 given  by (\ref{eq6}), and denote by $\alpha_\lambda>0$ the maximizer in $I$, we have: 
 $$\lambda x_\lambda=(1-\lambda) \alpha_\lambda (y_\lambda-x_\lambda)-(1-\lambda) \alpha_\lambda^2 x_\lambda.$$
 
The fact that  $I\subset J$ gives an advantage to player 2, which implies that no limit point of $x_\lambda$ and $y_\lambda$ can be greater than 1/2: 

\begin{lem}\label{lem5} Let $\lambda_n$ be a vanishing sequence of discount factors, such  that $y_{\lambda_{n}}$ and $x_{\lambda_{n}}$ converge to $v$ in $[0,1]$. Then $v\leq 1/2,$  $y_{\lambda_n}-x_{\lambda_n} \sim 2 \sqrt{\lambda_n} (1-v)$  and $\beta_{\lambda_n}\sim \sqrt{\lambda_n}.$
\end{lem}

\noindent {\bf Proof:} Using (\ref{eq8}) and the definition of $\alpha_\lambda$,   we get $x_\lambda(\lambda+\alpha_\lambda^2)= \lambda x_\lambda \alpha_\lambda^2+2 \alpha_\lambda  \sqrt{\lambda} \sqrt{1-\lambda}(1-y_\lambda).$
 Since $\lambda+\alpha_\lambda^2- 2 \alpha_\lambda \sqrt{\lambda}\geq 0$, we obtain:
 $ \lambda x_\lambda \alpha_\lambda^2+2 \alpha_\lambda  \sqrt{\lambda} \sqrt{1-\lambda}(1-y_\lambda)\geq 2 \alpha_\lambda \sqrt{\lambda} x_\lambda.$ Dividing by $\alpha_\lambda \sqrt{\lambda}$ and passing to the limit gives $v\leq 1/2$. The asymptotics are given by (\ref{eq8}) and the definition of $\beta_{\lambda_n}.$

\begin{lem}\label{lem6} 
 Let $\lambda_n$ be a vanishing sequence of discount factors such that $\sqrt{\lambda_n}\in  I$ for each $n$. Then  $y_{\lambda_n}$ and $x_{\lambda_n}$ converge to 1/2.
  \end{lem}
\noindent{\bf Proof:} By considering a converging subsequence we can assume that $y_{\lambda_n}$ and $x_{\lambda_n}$ converge to some $v$ in $[0,1]$. By the previous lemma, $v\leq 1/2$, and we have to show that $v\geq 1/2$. We have for each $\lambda$ in the subsequence, since player 1 can choose to play $\alpha=\sqrt{\lambda}$:
$$\lambda x_\lambda\geq (1-\lambda) \sqrt{\lambda}(y_\lambda-x_\lambda) -(1-\lambda) \lambda x_\lambda.$$
 Dividing by $\lambda$ and passing to the limit, we obtain   $v\geq 1/2$.
 
 \begin{lem} \label{lem7} Let $\lambda_n$ be a vanishing sequence of discount factors such that for each $n$,  the open interval $(\frac{1}{2}\sqrt{\lambda_n}, 2 \sqrt{\lambda_n})$ does not intersect $I$. Then    $\limsup_n y_{\lambda_n}\leq 4/9$. \end{lem}

\noindent{\bf Proof:} Suppose that (up to a subsequence) $x_{\lambda_n}$ and $y_{\lambda_n}$ converges to some $v\geq 4/9$. It is enough to show that $v=4/9$. We know that $v\leq 1/2$ by lemma \ref{lem5}. Consider again  the maximization problem of player 1 given  by (\ref{eq6}), and denote by $\alpha^*(\lambda)=\frac{y_\lambda-x_\lambda}{2(x_\lambda-z_\lambda)}>0$ the argmax of the unconstrained problem if player 1 could choose any $\alpha\geq 0$.  Since $\alpha^*(\lambda)\sim \sqrt{\lambda}\frac{1-v}{v}$ we have  $ \frac{1}{2}\sqrt{\lambda}\leq \alpha^*(\lambda)\geq 2\sqrt{\lambda}$ for $\lambda$ small in the sequence. By assumption $(\frac{1}{2}\sqrt{\lambda}, 2 \sqrt{\lambda})$ contains no point in $I$ and the objective function of player 1 is increasing from 0    to $\alpha^*(\lambda)$ and decreasing after $\alpha^*(\lambda)$. There are 2 possible cases:

If  $\alpha_\lambda\leq \frac{1}{2}\sqrt{\lambda}$ we have 
 $\lambda x_\lambda \leq \frac{1}{2} (1-\lambda) \sqrt{\lambda}(y_\lambda-x_\lambda) - \frac{1}{4}(1-\lambda)\lambda  -x_\lambda.$
Dividing by $\lambda$ and passing to the limit gives: $v\leq 1-v-\frac{1}{4}v$, i.e. $v\leq \frac{4}{9}.$

Otherwise, $\alpha_\lambda> 2 \sqrt{\lambda}$ and we have
 $\lambda x_\lambda \leq 2 (1-\lambda) \sqrt{\lambda}(y_\lambda-x_\lambda) -4(1-\lambda)\lambda  x_\lambda.$
Again, dividing by $\lambda$ and passing to the limit gives: $v\leq 4(1-v)- {4}v$, i.e. $v\leq \frac{4}{9}.$\\

Finally, the sequence $\lambda_n=\frac{1}{2^{2n}}$ satisfies the assumption of lemma \ref{lem6}, and the sequence $\lambda_n=\frac{1}{2^{2n+1}}$ satisfies those of lemma \ref{lem7}. This is enough to conclude that $v_\lambda(k_0)$ and $v_\lambda(k_1)$ have no limit when $\lambda$ goes to 0.

\subsection{The CV of $(v_n)_n$ and $(v_\lambda)_\lambda$ are equivalent.}

The equivalence between the uniform convergence of  $(v_n)_n$ and  $(v_\lambda)_\lambda$, which holds in general in 1-player games, has been recently proved (Ziliotto 2015) to extend to a large class of stochastic games.

It applies in particular to the following setup. Assume the set of states $K$ and the set of actions $I$ and $J$ are compact metric spaces, that the transition $q: K \times I \times J \longrightarrow \Delta(K)$ and the payoff $g: K \times I \times J \longrightarrow \R$ are jointly continuous. Together with an initial state $k$
, $(K,I,J,q,g)$ define a stochastic game. Then one can show that for each $n$ and each $\lambda$ the value of the $n$-stage game $v_n(k)$ and $v_\lambda(k)$ exist and satisfy the Shapley equations: $\forall n \geq 0, \forall \lambda \in (0,1],\; \forall k \in K, $

\begin{eqnarray*} 
 (n+1)\;  v_{n+1} (k)& =&\sup_{x \in \Delta(I)} \inf_{y \in \Delta(J)}  \left(g(k,x,y))+    n \; \E_{q(k,x,y)}  (v_n)\right),\\
& =&\inf_{y \in \Delta(J)}  \sup_{x \in \Delta(I)}  \left(g(k,x,y))+    n \; \E_{q(k,x,y)}  (v_n)\right).
\end{eqnarray*}

 \begin{eqnarray*} 
 v_{\lambda} (k)& =& \sup_{x \in \Delta(I)} \inf_{y \in \Delta(J)}  \left(\lambda \; g(k,x,y)+(1-\lambda)  \; \E_{q(k,x,y)}  (v_{\lambda}) \right), \\
 & =&  \inf_{y \in \Delta(J)}\sup_{x \in \Delta(I)}  \left(\lambda \; g(k,x,y)+(1-\lambda)  \; \E_{q(k,x,y)}  (v_{\lambda}) \right),
 \end{eqnarray*}

  \end{document}